\documentclass{article}

\usepackage[all]{xy}
 \usepackage[usenames,dvipsnames]{pstricks}
 \usepackage{epsfig}
 \usepackage{pst-grad} 
 \usepackage{pst-plot} 
\usepackage{ebezier}

\usepackage[
        colorlinks, citecolor=red,
        backref,
        pdfauthor={MDJ,JK},
        pdftitle={S1}
]{hyperref}

\usepackage{latexsym}\usepackage{amsfonts}\usepackage{amssymb}\usepackage{amsthm}
\usepackage{amsmath}\usepackage{newlfont}\usepackage{graphicx}\usepackage{doc}
\usepackage{color}\usepackage{multicol}\usepackage{epic}\usepackage{eepic}
\usepackage{ebezier}
\newtheorem{thrm}{Theorem}
\newtheorem{lem}[thrm]{Lemma}
\newtheorem{cor}[thrm]{Corollary}

\newtheorem{rem}[thrm]{Remark}
\newtheorem{prop}[thrm]{Proposition}

\newtheorem*{namedthm}{\namedthmname}
\newcounter{namedthm}

\makeatletter
\newenvironment{named}[1]
  {\def\namedthmname{#1}%
   \refstepcounter{namedthm}%
   \namedthm\def\@currentlabel{#1}}
  {\endnamedthm}
\makeatother

\include{epsf}

\def \Dj{\mbox{\raise0.3ex\hbox{-}\kern-0.4em D}}

\begin{document}

\title{Polynomial entropy of induced maps of circle and interval homeomorphisms}

\author{Ma\v{s}a \Dj ori\'c\\
Matemati\v{c}ki institut SANU\\
Knez Mihailova 36\\
11000 Beograd\\
Serbia\\
masha@mi.sanu.ac.rs \and
Jelena Kati\'c\thanks{Corresponding author: Jelena Kati\'c.} \\
Matemati\v cki fakultet\\
Studentski trg 16\\
11000 Beograd\\
Serbia\\
jelenak@matf.bg.ac.rs
}

\maketitle

\begin{abstract} We compute the polynomial entropy of the induced maps on hyperspace for a homeomorphism $f$ of an interval or a circle with finitely many non-wandering points.
\end{abstract}

\medskip

\section{Introduction}

Every continuous map on a compact metric space $X$ induces a continuous map (called the induced map) on the hyperspace $2^X$ of all closed subsets. If $X$ is connected, i.e.\ continuum, we consider the hyperspace $C(X)$ of subcontinua of $X$, (which is also a continuum). One can also consider the hyperspace $X^{*k}$ of all nonempty subsets with at most $k$ points (for $k\in\mathbb{N}$). A natural question is to find some relations between the given (individual) dynamics on $X$ and the induced one (collective dynamics) on the hyperspace. Various results in this direction were obtained in the last decades. Without attempting to give complete references, we mention just a few: Borsuk and Ulam~\cite{BU}, Bauer and Sigmund~\cite{BS}, Rom\'an-Flores~\cite{RF}, Banks~\cite{B}, Acosta, Illanes and M\'endez-Lango~\cite{AIM}.

The topological entropy of the induced map was studied by Kwietniak and Oprocha in~\cite{KO}, Lampart and Raith~\cite{LR}, Hern\'andez and M\'endez~\cite{HM}, Arbieto and Bohorquez~\cite{AB} and others.

In~\cite{LR} the authors showed that, if $f$ is an interval or a circle homeomorphism, the topological entropy of the induced map on the hyperspace of subcontinua is zero. This is also obtained as a corollary in~\cite{AB}, for Morse-Smale diffeomorphisms of a circle.

One of the measures of complexity of a system with zero topological entropy is the polynomial entropy. The notion of the polynomial entropy was first introduced by Marco in~\cite{M1} and~\cite{M2} in the context of Hamiltonian integrable systems. It was further investigated in different contexts by Labrousse~\cite{L1,L2,L3}, Labrousse and Marco~\cite{LM}, Bernard and Labousse~\cite{BL}, Artigue, Carrasco--Olivera and Monteverde~\cite{AOM}, Haseaux and Le Roux~\cite{HL}, Roth, Roth and Snoha~\cite{RRS}, Correa and de Paula~\cite{CP} etc.

As opposed to the topological entropy, which depends only on the dynamics restricted to the non-wandering set, the wandering set is visible to the polynomial entropy. In the case when non-wandering set is finite (for example Morse-Smale systems), there is a technique for computing the polynomial entropy developed in~\cite{HL} by Hauseux and Le Roux. Originally, they invented a simple coding procedure for homeomorphisms with only one non-wandering point, where the polynomial entropy is particularly well adapted, since the growth of the number of wandering orbits is always at least linear and at most polynomial. Hauseux and Le Roux also proved that the polynomial entropy localizes near a certain finite set (singular set), in order to compute the polynomial entropy of Brouwer homeomorphisms. This method was slightly generalized in~\cite{KP}, to the case of a map which is only continuous and the non-wandering set is finite. The coding procedure was also used in~\cite{CP} for the computation of the polynomial entropy of Morse-Smale systems on surfaces.

In this paper we compute the polynomial entropy of the induced maps $C(f)$, $f^{*k}$ and $2^f$ for a homeomorphism $f$ of a circle or an interval with a finite non-wandering set. The polynomial entropy of such a homeomorphism of an interval is known to be $1$, and of a circle is ether $0$ or $1$ (see~\cite{L3}). The hyperspace of subcontinua of an one-dimensional space is quite simple and can be identified with a two-dimensional object. Our computation uses the coding method and reduction to the singular sets.

Denote by $h_{pol}(\cdot)$ the polynomial entropy of a map and by $C(f)$, $f^{*k}$ and $2^f$ the induced maps on $C(X)$, $X^{*k}$ and $2^X$, respectively. These are the statements of our results.

\begin{named}{Theorem A}\label{lem:A}
{\it Let $f:[0,1]\to [0,1]$ be a homeomorphism with a finite non-wandering set. Then $h_{pol}(C(f))=2$, $h_{pol}(f^{*k})=k$ and $h_{pol}(2^{f})=\infty$.}\qed
\end{named}

\medskip

\begin{named}{Theorem B}\label{lem:B}
{\it Let $f:S^1\to S^1$ be a homeomorphism with a finite non-wandering set. Then $h_{pol}(C(f))=2$, $h_{pol}(f^{*k})=k$ and $h_{pol}(2^{f})=\infty$.}\qed
\end{named}

\section{Preliminaries}

\subsection{Hyperspaces and induced maps} For a compact metric space $(X,d)$, the hyperspace $2^X$ is the set of all nonempty closed subsets of $X$. The topology on $2^X$ is induced by the Hausdroff metric
$$d_H(A,B):=\inf\{\varepsilon\mid A\subset U(B,\varepsilon),\;B\subset U(A,\varepsilon)\},$$
where
$$U(A,\varepsilon):=\{x\in X\mid d(x,A)<\varepsilon\}.$$
The space $2^X$ is compact and the topology induced by $d_H$ is the same as the Vietoris topology.

We will also consider two closed subspaces of $2^X$, with the induced metric. The first one is $X^{*k}$, the space of all finite subsets of cardinality at most $k$, with the same topology. The set $X^{*k}$ is called $k$-fold symmetric product of $X$.

If $X$ is also connected (i.e.\ a continuum), then the set $C(X)$ of all connected and closed nonempty subsets of $X$ is also a continuum. The set $C(X)$ is called the hyperspace of subcontinua of $X$.

If $f:X\to X$ is continuous, then it induces continuous maps
$$\begin{aligned}&2^f:2^X\to 2^X,\quad &2^f(A):=\{f(x)\mid x\in A\}\\
&C(f):C(X)\to C(X),\quad &C(f)(A):=\{f(x)\mid x\in A\}\\
&f^{*k}:X^{*k}\to X^{*k},\quad &f^{*k}(A):=\{f(x)\mid x\in A\}.
\end{aligned}$$
If $f$ is a homeomorphism, so are $2^f$, $C(f)$ and $f^{*k}$.

\subsection{Polynomial entropy and coding}\label{subsec:entropy}

Suppose that $X$ is a compact metric space, and $f:X\rightarrow X$ is continuos. Denote by $d_n^f(x,y)$ the dynamic metric (induced by $f$ and $d$):
$$
d_n^f(x,y)=\max\limits_{0\leq k\leq n-1}d(f^k(x),f^k(y)).
$$
Fix $Y\subseteq X$. For $\varepsilon>0$, we say that a finite set $E\subset X$ is $(n,\varepsilon)$-separated if for every $x,y \in E$ it holds $d_n^f(x,y)\geq \varepsilon$. Let $S(n,\varepsilon;Y)$ denotes the maximal cardinality of an $(n,\varepsilon)$-separated set $E$, contained in $Y$.

\begin{def}\label{def:pol_ent} The polynomial entropy of the map $f$ on the set $Y$ is defined by
$$
h_{pol}(f;Y)=\lim\limits_{\varepsilon \rightarrow 0}\limsup\limits_{n\rightarrow \infty}\frac{\log S(n,\varepsilon;Y)}{\log n}.
$$
\end{def}
If $X=Y$ we abbreviate $h_{pol}(f):=h_{pol}(f;X)$.
The polynomial entropy, as well as the topological entropy, can also be defined via coverings with sets of $d^n_f$-diameters less than $\varepsilon$, or via coverings by balls of $d_n^f$-radius less than $\varepsilon$, see \cite{M1}. We list some properties of the polynomial entropy that are important for our computations:
\begin{itemize}
\item $h_{pol}(f^k)=h_{pol}(f)$, for $k\ge 1$
\item if $Y\subset X$ is a closed, $f$-invariant set, then $h_{pol}(f;Y)=h_{pol}(f|_Y)$
\item if $Y=\bigcup_{j=1}^mY_j$ where $Y_j$ are $f$-invariant, then $h_{pol}(f;Y)=\max\{h_{pol}(f;Y_j)\mid j=1,\ldots,m\}$
\item If $f:X\to X$, $g:Y\to Y$ and $f\times g:X\times Y\to X\times Y$ is defined as $f\times g (x,y):=(f(x),g(y))$, then $h_{pol}(f\times g)=h_{pol}(f)+h_{pol}(g)$
\item $h_{pol}(f)$ does not depend on a metric but only on the induced topology
\item $h_{pol}(\cdot)$ is a conjugacy invariant (meaning if $f:X\to X$, $g:X'\to X'$, $\varphi:X\to X'$ is a homeomorphism of compact spaces with and $g\circ\varphi=\varphi\circ f$, then $h_{pol}(f)=h_{pol}(g)$).
\end{itemize}

A point $p\in X$ is wandering if there exists a neighbourhood $U\ni p$ such that $f^n(U)\cap U=\emptyset$, for $n\ge 1$.\label{wan-pt}

A point that is not wandering is said to be non-wandering. We denote the set of all non-wandering points by $NW(f)$. The set $NW(f)$ is  closed and $f$-invariant.

We now give a brief description of the computation of the polynomial entropy for maps with a finite non-wandering set, by means of a coding and a local polynomial entropy. This construction was first done in~\cite{HL} for homeomorphisms with only one non-wandering (hence fixed) point, and then modified in~\cite{KP} for continuous maps with finitely many non-wandering points. Let $Y$ be any $f$-invariant subset of $X$.

We first define a coding relative to a family of sets $\mathcal{F}$.  Let
$$\mathcal{F}=\{Y_1,Y_2,\ldots,Y_L\}$$ where $Y_j\subseteq X\setminus\mathrm{NW}(f)$ and
$$Y_{\infty}:=Y\setminus\bigcup_{j=1}^LY_j.$$ Let $\underline{x}=(x_0,\ldots,x_{n-1})$ be a finite sequence of elements in $Y$. We say that a finite sequence $\underline{w}=(w_0,\ldots,w_{n-1})$ of elements in $\mathcal{F}\cup\{Y_\infty\}$ is a {\it coding} of $\underline{x}$ relative to $\mathcal{F}$ if $x_j\in w_j$. We will refer to $\underline{w}$ as a word and to $w_j$ as a letter.

Let $\mathcal{A}_n(\mathcal{F};Y)$ be the set of all codings of all orbits
$$(x,f(x),\ldots,f^{n-1}(x))$$ of length $n$ relative to $\mathcal{F}$, for all $x\in Y$. If
$\sharp \mathcal{A}_n(\mathcal{F};Y)$ denotes the cardinality of $\mathcal{A}_n(\mathcal{F};Y)$, we define the {\it polynomial entropy of $f$, on the set $Y$, relative to the family} $\mathcal{F}$ as the number:
$$h_{pol}(f,\mathcal{F};Y):=\limsup_{n\to\infty}\frac{\log \sharp \mathcal{A}_n(\mathcal{F};Y)}{\log n}.$$

We abbreviate $h_{pol}(\mathcal{F};Y):= h_{pol}(f,\mathcal{F};Y)$ whenever there is no risk of confusion.

For $Z\subseteq X\setminus\operatorname{NW}(f)$, set
$$M(Z):=\sup_{x\in X}\sharp\{n\mid f^n(x)\in Z\}.$$
If $Z\subseteq X\setminus\operatorname{NW}(f)$ is compact, the number $M(Z)$ is finite, since $Z$ can be covered with a finite number of open wandering sets, and every orbit can intersect a wandering set at most once.

We will use the following property of $h_{pol}(\mathcal{F};Y)$ in order to localize our compuutation to a singular set.

\begin{prop}\label{prop:monotonicity}[Proposition 3.2 in~\cite{KP}] Let $\mathcal{F}$ and $\mathcal{F}'$ be two families of subsets of $X\setminus\mathrm{NW}(f)$ with $M(\cup\mathcal{F})<+\infty$. Let $Y\subset X$ be an $f$-invariant subset with exactly one non-wandering point. If for every $Y_j'\in\mathcal{F}'$ there exists $Y_j\in\mathcal{F}$ such that $Y'_j\subseteq Y_j$, then $h_{pol}(\mathcal{F}';Y)\le h_{pol}(\mathcal{F};Y).$\qed
\end{prop}

Next we define the local polynomial entropy for a finite set
$$\mathcal{S}=\left\{x_1,x_2,\dots,x_l\right\}\subset X\setminus NW(f).$$
We choose a decreasing sequence of neighbourhoods $\{U_{j,n}\}_{n \in \mathbb{N}}$ of $x_j\in \mathcal{S}$ which form a basis of neighbourhoods of $x_j$. It follows from Proposition~\ref{prop:monotonicity} that the sequence
$$\left\{h_{pol}\left(\{U_{1,n},U_{2,n},\dots, U_{l,n}\};Y\right)\right\}_{n \ge 1}$$ is decreasing and converges, as well as that its limit does not depend on the choice of neighbourhoods. Define:
$$
h_{pol}^{loc}(\mathcal{S};Y):=\lim\limits_{n\rightarrow \infty}h_{pol}\left(\{U_{1,n},U_{2,n},\dots, U_{l,n}\};Y\right).
$$

Finally, we relate the polynomial entropy to a singular set.

We say that the subsets $U_1,\ldots,U_L$ of $X\setminus\mathrm{NW}(f)$ are {\it mutually singular} if for every $M>0$, there exist $x$ and positive integers $n_1,\ldots,n_L$ such that
$$f^{n_j}(x)\in U_j,\quad |n_i-n_j|>M, \quad\text{for}\;i\neq j.$$
The points $x_1,\ldots,x_L$ are mutually singular if every family of respective neighbourhoods $U_1,\ldots,U_L$, $U_j\ni x_j$ is mutually singular. We say that a finite set is singular if it consists of mutually singular points.

\begin{prop}\label{prop:sing}{\rm [Propositions 3.2 and 3.3 in~\cite{KP}]} Let $Y\subseteq X$ be an $f$-invariant set containing exactly one non-wandering point. Then it holds:
\begin{itemize}
\item[(a)] $h_{pol}(f;Y)=\sup\{h_{pol}(\{K\};Y)\mid {K\subseteq X\setminus\operatorname{NW}(f),\,K \,\text{compact}}\}$
\item[(b)] $h_{pol}(f;Y)=\sup\{h_{pol}(\{K\cap Y\};Y)\mid{K\subseteq X\setminus\operatorname{NW}(f),\,K \,\text{compact}}\}$
\item[(c)] $h_{pol}(f;Y)=
\sup\big\{h_{pol}^{loc}(\mathcal{S};Y)\mid\mathcal{S}\subset\overline{Y},\;\mathcal{S} \,\text{singular}\big\}$.\qed
\end{itemize}
\end{prop}

\medskip

The following corollary is of importance for our result, so we will prove it.

\begin{cor}\label{cor:upper}
The polynomial entropy $h_{pol}(f;Y)$ is bounded from above by the maximal cardinality of a singular set contained in $\overline{Y}$.
\end{cor}
\noindent{\it Proof.} The local polynomial entropy of a finite set $\{x_1,\ldots,x_l\}\subset{X}\setminus\operatorname{NW}(f)$ is bounded from above by its cardinality. Indeed, one can choose wandering neighbourhoods $Y_j\ni x_j$, such that every letter $Y_j$ appears in the coding of any orbit $(x,f(x),\ldots,f^{n-1}(x))$ at most once. Therefore there are at most $n(n-1)\cdots(n-l+1)\le n^l$ possible codings. We conclude
$$h_{pol}^{loc}(\{x_1,\ldots,x_l\};Y)\le h_{pol}(\{Y_1,\ldots,Y_l\};Y)=\limsup_{n\to\infty}\frac{\log\;\sharp\;\mathcal{A}_n(\{Y_1,\ldots,Y_l\};Y)}{\log n}\le l.$$\qed

In particular, if the maximal cardinality of a singular set is two, we have a more precise statement.

\begin{prop}\label{prop:eq2} Let $Y$ be as above. Suppose that the maximal cardinality of a singular set contained in $\overline{Y}$ equals 2. If for every singular set $\mathcal{S}=\{x_1,x_2\}\subset Y$ and every open $Y_j\ni x_j$, $Y_j\subset X\setminus NW(f)$, there exists a positive integer $L$ such that for all $n\ge L$ it holds $f^n(Y_1)\cap Y_2\neq\emptyset$, then $h_{pol}(f;Y)=2$.
\end{prop}
\noindent{\it Proof.} It follows from Corollary~\ref{cor:upper} that $h_{pol}(f;Y)\le 2$. To prove the other inequality, note that the assumed condition on $Y_1$ and $Y_2$ implies that for all positive integer $m$ and for all $n\ge L$, there exists $x$ with $f^m(x)\in Y_1$, $f^{m+n}(x)\in Y_2$. Therefore, for every $n\ge L$ there exists a coding of a form
$$(\underbrace{Y_\infty,\ldots,Y_\infty}_{m},Y_1,\underbrace{Y_\infty,\ldots,Y_\infty}_{n},Y_2,Y_\infty,\ldots,Y_\infty)$$ where $m$ is any number and $n\ge L$. We get
$$\sharp\;\mathcal{A}_n(\{Y_1,Y_2\};Y)\ge\sum_{m=0}^{n-L-1}(n-L-m)\sim \frac12n^2,$$
so $h_{pol}(\{Y_1,Y_2\};Y)\ge 2$. Since this holds for all $Y_j$ with $x_j\in Y_j\subset X\setminus NW(f)$, the same is true for $h_{pol}(\{x_1,x_2\};Y)$.\qed

\smallskip

We will use the following notations:
$$\mathcal{O}(x):=\{f^n(x)\mid n\in\mathbb{Z}\}$$ for an orbit of a point $x$ and
$$W^s(p):=\{x\in X\mid f^n(x)\to p,n\to\infty\}$$ for the stable manifold of a (fixed) point $p$.

\section{Proofs}

In order to use the methods described in Subsection~\ref{subsec:entropy} we need to establish that the sets $NW(C(f))$ and $NW\left(f^{*k}\right)$ are finite. Since $h_{pol}(f^m)=h_{pol}(f)$, we can assume that $NW(f)=Fix(f)$.

\begin{prop}
Let $f:[0,1]\to[0,1]$ or $f:S^1\to S^1$ be a homeomorphism such that the set $NW(f)=Fix(f)$ is finite. Then $NW(C(f))=Fix(C(f))$ and $NW\left(f^{*k}\right)=Fix\left(f^{*k}\right)$ are also finite for every $k\ge 1$.
\end{prop}

\noindent{\it Proof.} Let us prove the proposition for $C(f)$, $f:I\to I$. The other case can be proved analogously.
The set $Fix(C(f))$ is finite as it consists of all intervals $[a,b]\subseteq I$, where $a,b\in Fix(f)$. Note that the condition $NW(f)=Fix(f)$ implies that $f$ is increasing. Suppose that there exists $[x,y]\in NW(C(f))\setminus Fix(C(f))$ and that $x\notin Fix(f)$. Since $NW(f)=Fix(f)$, there exist $a,b\in Fix(f)$ (possibly equal) with
$$f^n(x)\to a, \quad f^n(y)\to b, \quad n\to\infty.$$ Let $\varepsilon<d(x,a)/2$. For $n\ge n_0$ it holds $d(f^n(x),a)<\varepsilon$, $d(f^n(y),b)<\varepsilon$ so we have
$$\begin{aligned}
d_H(C(f)^n([x,y]),[x,y])&\ge
d_H([a,b],[x,y])-d_H\left(C(f)^n([x,y]),[a,b]\right)\\
&\ge d(a,x)-\varepsilon>\varepsilon.
\end{aligned}$$
\qed

\subsection{Proof of \ref{lem:A}}

Let $I:=[0,1]$. We will first compute the polynomial entropy of $C(f)$. We can identify the space $C(I)$ with the set
$$\{(x,y)\in[0,1]^2\mid 0\le x\le y\le 1\},$$ which is the upper triangle $A$ in the square $[0,1]\times [0,1]$. The homeomorphism $\varphi:C(I)\to A$ is given by $\varphi:[x,y]\to(x,y)$. The map $C(f)$ is conjugated to $f\times f|_A$ via $\varphi$ (where $f\times f:[0,1]\times [0,1]\to[0,1]\times [0,1]$, $f(x,y)=(f(x),f(y))$).

Denote the lower triangle in $[0,1]\times [0,1]$ by $B$. Since both $A$ and $B$ are closed and $f\times f$-invariant (if not, they are $(f\times f)^2$-invariant), we have
$$h_{pol}(f\times f)=\max\{h_{pol}(f\times f|_A),h_{pol}(f\times f|_B)\}.$$ On the other hand
$$h_{pol}(f\times f)=2h_{pol}(f).$$
Since $f\times f|_A$ and $f\times f|_B$ are topologically equivalent systems (the map $(x,y)\mapsto(y,x)$ realizes a conjugacy), we have
$$h_{pol}(C(f))=2h_{pol}(f).$$

It is not hard to see that $h_{pol}(C(f))=2$, because $h_{pol}(f)=1$. Indeed, suppose that $f$ is increasing (if not, $f^2$ is). Let $0=p_0<p_1<\ldots<p_k=1$ denote the fixed points. We apply again $h_{pol}(f)=\max_j h_{pol}\left(f|_{[p_{j-1},p_j]}\right)$.
Since $f|_{[p_{j-1},p_j]}$ possesses a non-wandering point, it holds $h_{pol}(f|_{[p_{j-1},p_j]})\ge 1$ (see, for example~\cite{L3}). It is easy to see that $f|_{[p_{j-1},p_j]}$ satisfies all the conditions in Proposition~\ref{prop:sing}, for $Y=(p_{j-1},p_j]$, therefore also in Corollary~\ref{cor:upper}, as well as that $f$ has only one singular point in $(p_{j-1},p_j]$. So we obtain
$$h_{pol}\left(f|_{[p_{j-1},p_j]}\right)=\max\left\{h_{pol}\left(f|_{[p_{j-1},p_j]};(p_{j-1},p_j]\right),h_{pol}\left(f|_{\{p_{j-1}\}}\right)\right\}\le 1.$$ (For a different, more explicit proof of $h_{pol}(f|_{[p_{j-1},p_j]})=1$ see~\cite{L3}.)

\medskip

Let us prove that $h_{pol}\left(f^{*k}\right)=k$. Since $f^{\times k}:=f\times\ldots\times f$ and $f^{*k}$ are semi-conjugated via:
$$\pi:I^k\to I^{*k},\quad \pi:(x_1,\ldots,x_k)\mapsto\{x_1,\ldots,x_k\},$$ we have
\begin{equation}\label{eq:<=k}
h_{pol}\left(f^{*k}\right)\le h_{pol}\left(f^{\times k}\right)=k.
\end{equation}

We want to prove the other inequality. For a permutation $\sigma$ of $\{1,\ldots,k\}$, define
$$A_\sigma:=\{(x_1,\ldots,x_k)\in I^k\mid x_{\sigma(1)}\le\ldots\le x_{\sigma(k)}\}.$$ As before, we see that $I^{k}=\bigcup_\sigma A_\sigma$, $A_\sigma$'s are $f^{\times k}$-invariant, and all $f^{\times k}|_{A_\sigma}$ are mutually conjugated, so
$$h_{pol}\left(f^{\times k}\right)=\max\left\{h_{pol}\left(f^{\times k}|_{A_\sigma}\right)\right\}.$$ Therefore
$$h_{pol}\left(f^{\times k}|_{A_\sigma}\right)=k$$ for every $\sigma$.

Define:
$$\begin{aligned}A(k)&:=\{(x_1,\ldots,x_k)\in I^k\mid x_1\le\ldots\le x_k\}\\
\hat{A}(k,I)&:=\left\{\{x_1,\ldots,x_k\}\in I^{*k}\mid x_1<\ldots<x_k\right\}.\end{aligned}$$
Whenever there is no risk of confusion we will abbreviate $\hat{A}(k)=\hat{A}(k,I)$.
If $k=2$ we are done, since $\pi|_{A(2)}:A(2)\to I^{*2}$ is a homeomorphism of compact sets, so $h_{pol}\left(f^{*2}\right)=h_{pol}\left(f^{\times2}|_{A(2)}\right)=2$.
\medskip

For $k>2$ we need the following auxilirary fact.
\begin{lem}\label{lem} Suppose that $Fix(f)=\{0,1\}$. Then $h_{pol}\left(f^{*k};\hat{A}(k)\right)=k$.
\end{lem}

The rest of the proof follows easily from Lemma~\ref{lem}. Indeed, if $Fix(f)=\{0,1\}$ then we have:
$$h_{pol}\left(f^{*k}\right)\ge h_{pol}\left(f^{*k};{\hat{A}(k)})\right){=}k\;\stackrel{(\ref{eq:<=k})}\Rightarrow\;h_{pol}\left(f^{*k}\right)=k.$$ If $Fix(f)=\{0,p_1,\ldots,p_{k-1},1\}$ with $0<p_1<\ldots<p_{k-1}<1$, define
$$B:=\left\{\{x_1,\ldots,x_j\}\mid j\in\{1,\ldots,k\},\; x_i\in[0,p_1]\right\}\subset I^{*k}$$ and apply Lemma~\ref{lem} to $f|_{[0,p_1]}$. Since $B$ is $f^{*k}$-invariant, we conclude
$$h_{pol}\left(f^{*k}\right)\ge h_{pol}\left(f^{*k}|_B\right)=h_{pol}\left(\left(f|_{[0,p_1]}\right)^{*k}\right)\geqslant h_{pol}\left(\left(f|_{[0,p_1]}\right)^{*k};\hat{A}(k,[0,p_1])\right)=k.$$ From here and~(\ref{eq:<=k}) we obtain $h_{pol}\left(f^{*k}\right)=k$.\qed

\medskip

\noindent{\bf Proof of Lemma~\ref{lem}}. Consider the following (non-disjoint) partition of $A(k)$:
\begin{itemize}
\item $\tilde{A}(k):=\{(x_1,\ldots,x_k)\in I^k\mid x_1<\ldots<x_k\}$
\item $A_{j}(k):=\{(x_1,\ldots,x_k)\in I^k\mid x_1\le\ldots\le x_j=x_{j+1}\le\ldots\le x_k\}$, for $j=1,\ldots,k-1$.
\end{itemize}
It is obvious that $\tilde{A}(k)$ and $A_{j}(k)$ are $f^{\times k}$-invariant. Therefore
$$h_{pol}\left(f^{\times k}\right)=\max\left\{h_{pol}\left(f^{\times k};\tilde{A}(k)\right),h_{pol}\left(f^{\times k}|_{A_{j}(k)}\right)\right\}.$$
Notice that $f^{\times k}|_{A_{j}(k)}$ is conjugated to $f^{\times (k-1)}|_{A(k-1)}$, so
$$h_{pol}\left(f^{\times k}|_{A_{j}(k)}\right)=h_{pol}\left(f^{\times (k-1)}\right)=k-1.$$ Since $h_{pol}\left(f^{\times k}\right)=k$ we conclude that $h_{pol}\left(f^{\times k};\tilde{A}(k)\right)=k$.

Notice that
$$\pi|_{\tilde{A}(k)}:\tilde{A}(k)\to\hat{A}(k)$$ is a homeomorphism that establishes a conjugacy between $f^{\times k}|_{\tilde{A}(k)}$ and $f^{*k}|_{\hat{A}(k)}$. Although polynomial entropy is a conjugacy invariant only when the domain is compact (which the sets $\hat{A}(k)$ and $\tilde{A}(k)$ are not), we can indirectily prove that
$$
h_{pol}\left(f^{*k};{\hat{A}(k)})\right)=h_{pol}\left(f^{\times k};\tilde{A}(k)\right)=k.
$$
We wish to apply the coding method from Proposition~\ref{prop:sing}. Note that for $k>2$ the sets $\tilde{A}(k)$ and $\hat{A}(k)$ do not contain any non-wandering points. We can add the point $(0,\ldots,0)$ to the set $\tilde{A}(k)$ and $\{0\}$ to $\hat{A}(k)$, keeping the same notations: this will not change the entropy and $\pi$ will still be a homeomorphism between the two. In this way we achive that the assumptions from Proposition~\ref{prop:sing} are fulfilled.

Note that $\pi$ induces a bijection between the sets
$$\{K\cap\tilde{A}(k)\mid K\subset I^{k}\setminus NW\left(f^{\times k}\right),\;K\;\text{compact}\}$$ and
$$\{L\cap\hat{A}(k)\mid L\subset I^{*k}\setminus NW\left(f^{*k}\right),\;L\;\text{compact}\}.$$ If $\underline{w}=(w_0,\ldots,w_{n-1})$ is a coding of an orbit $\left(x,f^{\times k}(x).\ldots,(f^{\times k})^{n-1}(x)\right)$ in
$\tilde{A}(k)$ (consisting of letters $K\cap\tilde{A}(k)$ and $Y_\infty:=\tilde{A}(k)\setminus K$), then $(\pi(w_0),\ldots,\pi(w_{n-1}))$ is the coding of an orbit $\left(\pi(x),f^{*k}(\pi(x)),\ldots,(f^{*k})^{n-1}(\pi(x))\right)$ in $\hat{A}(k)$ (consisting of letters $\pi(K\cap\tilde{A}(k))$ and $Y_\infty:=\hat{A}(k)\setminus\pi(K\cap\tilde{A}(k))$), and vice versa. Therefore, for a fixed compact $K$, the sets $\mathcal{A}_n\left(\{K\cap\tilde{A}(k)\};\tilde{A}(k)\right)$ and $\mathcal{A}_n\left(\{\pi(K\cap\tilde{A}(k))\};\hat{A}(k)\right)$ have the same cardinality. Applying $(b)$ from Proposition~\ref{prop:sing} we finish the proof.

\medskip

Finally, to prove that $h_{pol}\left(2^f\right)=\infty$, we notice that $X^{*k}$ is a closed and $2^f$-invariant subset of $2^X$ and moreover $f^{*k}=2^f|_{X^{*k}}$ so
$$h_{pol}\left(2^f\right)\ge h_{pol}\left(f^{\times k}\right)=k$$ for every $k\in\mathbb{N}$.\qed

\subsection{Proof of \ref{lem:B}}

 Let us first compute the polynomial entropy of $C(f)$. For that reason, we will distinguish between the following possibilities:
\begin{itemize}
\item[(1)] the set $Fix(f)$ consists of only one point
\item[(2)] the set $Fix(f)$ has at least three different points
\item[(3)] the set $Fix(f)$ has exactly two points.
\end{itemize}

\begin{figure}
\centering
\includegraphics[width=11cm,height=5cm]{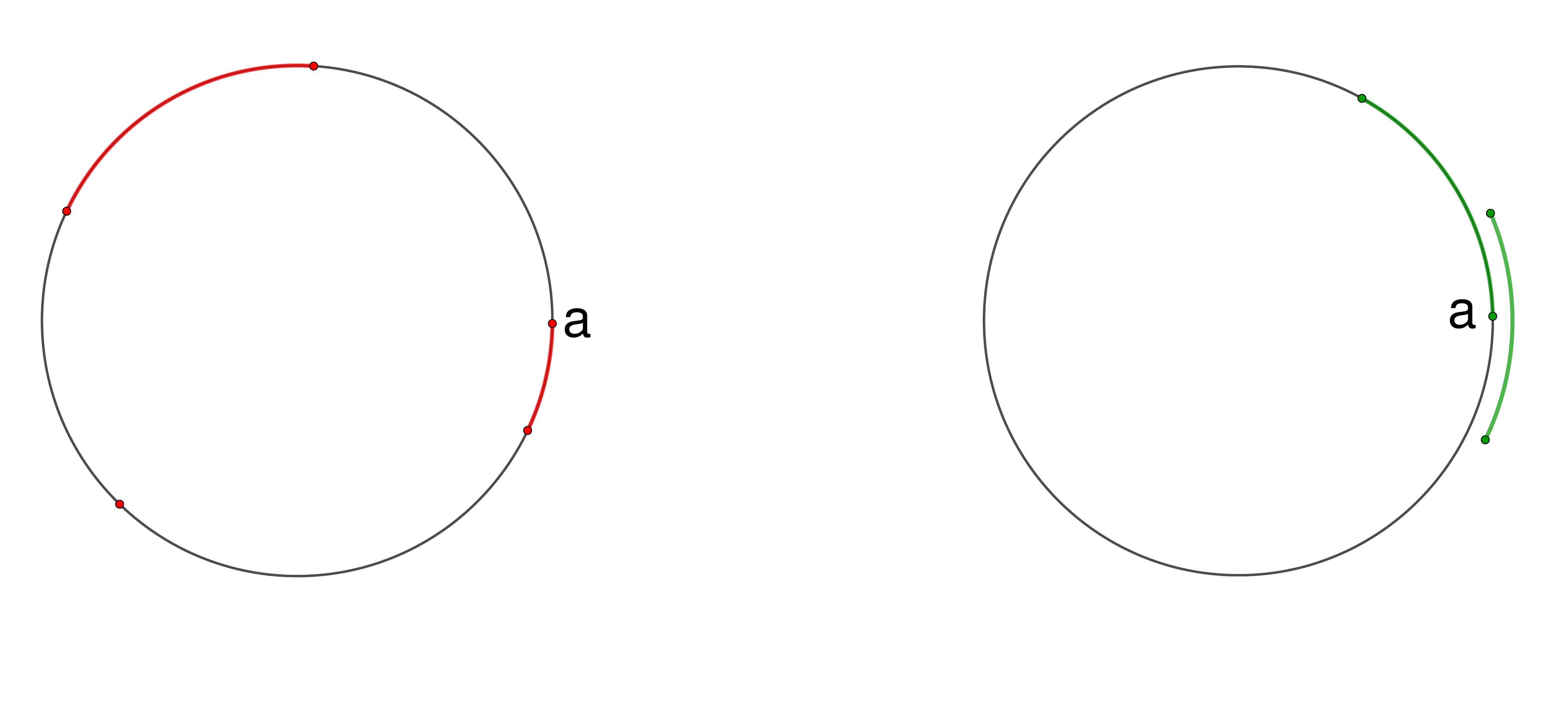}
\centering
\caption{Elements of set $A$ (on the left) and $B$ (on the right)}
\label{A&B}
\end{figure}

\noindent{\bf Case (1).} Since $f$ has only one fixed point, $f$ preserves the orientation of the circle. If $Fix(f)=\{a\}$, the continuum map $C(f)$ has only two non-wandering points - $\{a\}$ and $S^1$.
We will divide the set $C(S^1)$ into two closed invariant subsets:
\begin{itemize}
\item $A$ is the set of all $[x,y]$ such that the point $x$ is between points $a$ and $y$ counterclockwise, including degenerate cases when the two or all three points are equal (meaning that $[x,a]$, $\{x\}$, $[a,y]$ and $S^1$ are in $A$); notice that $[x,y]$ does not contain $a$ as an interior point,  for $x\neq y$
\item $B$ is the set of all $[x,y]$ such that the point $y$ is between points $a$ and $x$ counterclockwise, including degenerate cases when the two or all three points are equal (meaning that $[x,a]$, $\{a\}$, $[a,y]$ and $S^1$ are in $B$); notice that $[x,y]$ contains $a$
\end{itemize}
(see Figure~\ref{A&B}). In this way we know that, for $[x,y], [z,w]\in A$ or $[x,y], [z,w]\in B$ the following is true:
$$d(x,z)<r\quad\mbox{and}\quad d(y,w)<r \;\Rightarrow\; d_H([x,y],[z,w])<r.$$
(In general, this does not have to hold, since the arcs $[x,y]$ and $[z,w]$ may be "at the oposite sides" of $S^1$.)

\begin{figure}
\centering
\includegraphics[width=11cm,height=5cm]{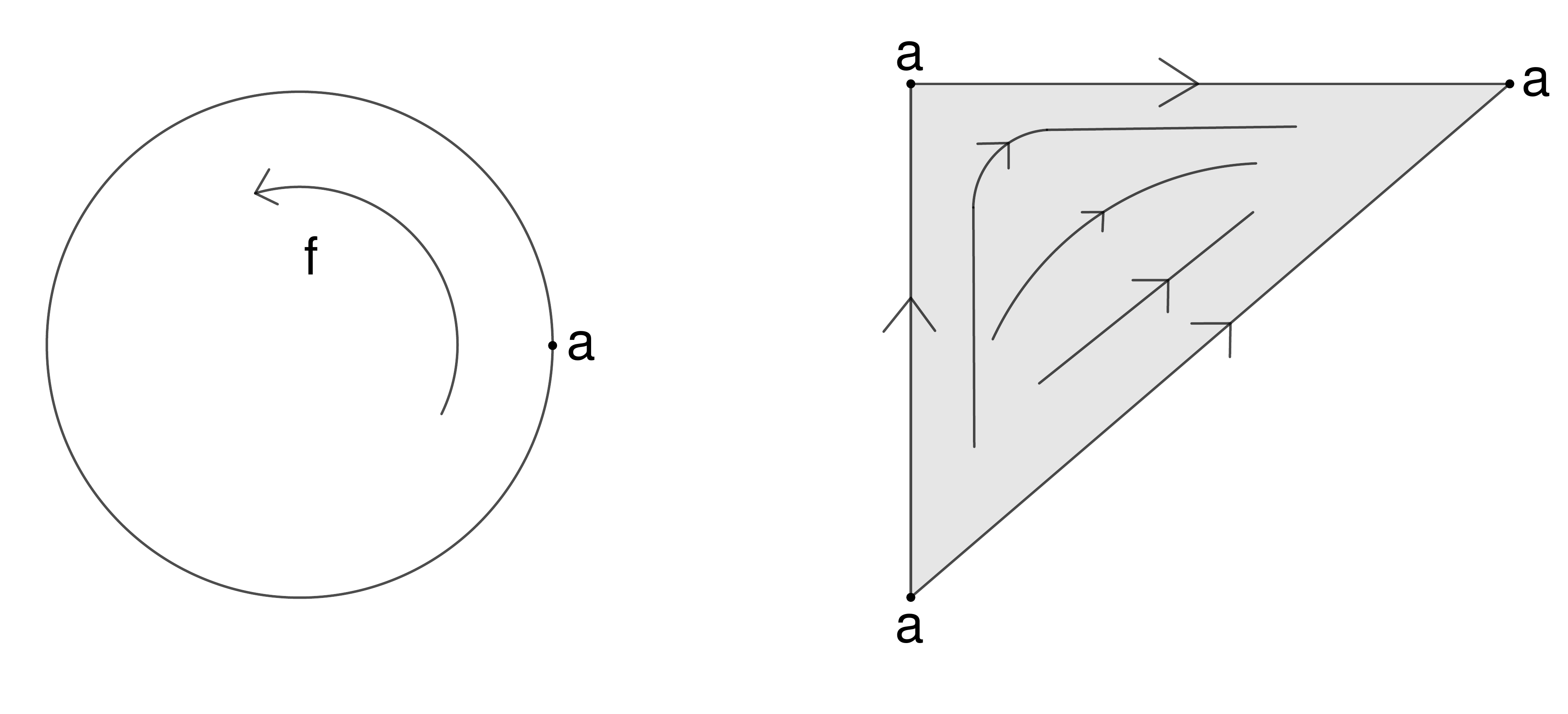}
\centering
\caption{Direction of $f$ (on the left) and dynamics of $C(f)|_A$ (on the right)}
\label{dir&dyn}
\end{figure}

 Suppose that $f$ moves the points in $S^1$ in the positive direction, as in Figure~\ref{dir&dyn}. The other case is treated in the same way.

 We will first consider the map $C(f)|_A$. Its dynamics is depicted in Figure~\ref{dir&dyn}. We have the following four possibilities:
 \begin{itemize}
 \item $x\neq a$ $\Rightarrow$ $C(f)^n([x,y])\stackrel{n\to\pm\infty}{\longrightarrow}\{a\}$
 \item $[x,y]=\{a\}$ $\Rightarrow$ $C(f)^n([x,y])\stackrel{n\to\pm\infty}{\longrightarrow}\{a\}$
 \item $x=a$, $y\neq a$ $\Rightarrow$ $C(f)^n([x,y])\stackrel{n\to\infty}{\longrightarrow}S^1$, $C(f)^n([x,y])\stackrel{n\to-\infty}{\longrightarrow}\{a\}$
 \item $[x,y]=S^1$ $\Rightarrow$ $C(f)^n([x,y])\stackrel{n\to\pm\infty}{\longrightarrow}S^1$.
 \end{itemize}

 We can divide $A$ into the sets $\{S^1\}$ and $Y:=A\setminus\{S^1\}$ and compute
 $$h_{pol}(C(f)|_A)=\max\left\{h_{pol}(C(f)|_A;Y),h_{pol}\left(C(f)|_A;\{S^1\}\right)\right\}=h_{pol}(C(f)|_A;Y)$$
since $ h_{pol}(C(f)|_A;\{S^1\})= h_{pol}(C(f)|_{\{S^1\}})=0$.

We claim that the arcs $[a,p]$, for $p\neq a$ and $[q,a]$ for $q\neq a$ are two mutually singular points $x_1$ and $x_2$ satisfying the conditions stated in Proposition~\ref{prop:eq2}.

Let us first prove that $[a,p]$ and $[q,a]$ are mutually singular. Fix an $\varepsilon>0$ $\varepsilon<\min\{d(a,p), d(a,q)\}$ and $M>0$.
Choose $y\in B(p,\varepsilon)\subset S^1$ arbitrarly. Since $f^n(y)\to a$, when $n\to\infty$, there exists a non-negative integer $n_1> M$ such that for $n\ge n_1$ it holds $d(f^n(y),a)<\varepsilon$.  Let $\alpha\in B(q,\varepsilon)\subset S^1$ be any point. Since $f^{-n}(\alpha)\to a$, when $n\to\infty$, there exsts $n_2\ge n_1$ with $d(f^{-n_2}(\alpha),a)<\varepsilon$. We can increase $n_2$ if necessary to obtain that the point $f^{-n_2}(\alpha)$ is between $a$ and $y$, and $f^{-n_2}(\alpha)\neq y$. Choose $x:=f^{-n_2}(\alpha)$.  Set $I:=[x,y]$.
We claim that the orbit of $I$ intersects the $\varepsilon$-balls around $[a,p]$ and $[q,a]$ in the times with difference greater than $M$. Indeed, we have:
$$d_H(I,[a,p])\le\max\{d(x,a),d(y,p)\}<\varepsilon,$$
so $I\in B([a,p],\varepsilon)$ and
$$d(f^{n_2}(I),[q,a])\le\max\{d(f^{n_2}(x),q),d(f^{n_2}(y),a\}=\max\{d(\alpha,q),d(f^{n_2}(y),a)\}<\varepsilon$$ therefore $f^{n_2}(I)\in B([q,a],\varepsilon)$.

The next step is to show that any two arcs except the ones of the form $[a,p]$ and $[q,a]$ cannot be mutually singular. Suppose that $[p_1,q_1]$ and $[p_2, q_2]$ are two different arcs such that $p_j\neq a$, $q_j\neq a$. If $[p_1,q_1]$ and $[p_2, q_2]$ are not on the same orbit of $C(f)$, it is enough to show that there exist neighbourhoods $U\ni[p_1,q_1]$ and $V\ni[p_2, q_2]$ such that
$$C(f)^n(U)\cap V=\emptyset,\quad \mbox{for all}\;n\in\mathbb{Z}.$$
Let $\varepsilon$ be a positive real number with the following properties:
\begin{itemize}
\item $\varepsilon< d_H\left([p_2,q_2],\mathcal{O}([p_1,q_1])\right)$, where $d_H\left([p_2,q_2],\mathcal{O}([p_1,q_1])\right)$
is the distance from $[p_2,q_2]$ to the orbit $\mathcal{O}([p_1,q_1])$, which is strictly positive, as $C(f)^n([p_1,q_1])$ converges to $\{a\}$, when $n\to\pm\infty$
\item the balls of radius $\varepsilon$ around $[p_2, q_2]$ and $\{a\}$ are disjoint.
\end{itemize}
Since $C(f)^n([p_1,q_1])$ converges to $\{a\}$, when $n\to\pm\infty$, and the same holds for any $[p,q]$ close enough to $[p_1,q_1]$, we can find $\delta>0$ and $n_0$ such that:
$$[p,q]\in B\left([p_1,q_1],\delta\right),\,|n|\ge n_0\quad\Rightarrow\quad C(f)^n\left([p,q]\right)\in B\left(\{a\},\varepsilon\right).$$ We can decrease $\delta$ if necessary to obtain:
$$[p,q]\in B\left([p_1,q_1],\delta\right),\;|n|\le n_0\quad \Rightarrow\quad C(f)^n\left([p,q]\right)\notin B\left([p_2,q_2],\varepsilon\right).$$ We conclude that the sets
$$U:=B\left([p_1,q_1]),\varepsilon\right),\quad V:=B\left([p_2,q_2],\delta\right)$$ have the desired properties, hence $[p_1,q_1]$ and $[p_2, q_2]$ are not mutually singular.

If $[p_1,q_1]$ and $[p_2, q_2]$ are on the same orbit, it is enough to show that there exist neighbourhoods $U\ni[p_1,q_1]$ and $V\ni[p_2, q_2]$, $U,V\subset Y\setminus NW(C(f))$ and $M>0$, such that:
\begin{equation}\label{eq:non-sing}
C(f)^{n}(U)\cap V\neq\emptyset\;\Rightarrow\; |n|<M,\quad \mbox{for all}\;n\in\mathbb{Z}.
\end{equation}
Indeed, it follows from~(\ref{eq:non-sing}) that $U$ and $V$ are not mutually singular, so neither are $[p_1,q_1]$ and $[p_2, q_2]$. So take $U\subset W^s(\{a\})\cap Y$, $V$ and $U'$ to be any three balls centered at $[p_1,q_1]$, $[p_2,q_2]$ and $\{a\}$ respectively, such that $V\cap U'=\emptyset$. There exists a non-negative integer $n_0$ such that, for all $n\in\mathbb{Z}$, $|n|\ge n_0$, $C(f)^n(U)\subset U'$. We see that for $|n|\ge n_0$ it holds $C(f)^n(U)\cap V=\emptyset$.

One can check in a similar way that for all other possibilities for $[p_1,q_1]$ and $[p_2, q_2]$ (except $[a,p]$, for $p\neq a$ and $[q,a]$ for $q\neq a$), the arcs $[p_1,q_1]$ and $[p_2, q_2]$ can not be mutually singular.

It remains to prove that for every two mutually singular points of the form $[a,p]$ and $[q,a]$, where $p,q\neq a$ and every open $Y_1\ni [a,p]$, $Y_2\ni[q,a]$, $Y_j\subset A\setminus NW(C(f))$, there exists a positive integer $L$ such that for all $n\ge L$ it holds $C(f)^n(Y_1)\cap Y_2\neq\emptyset$. Then we are able to apply Proposition~\ref{prop:eq2} and finish the proof.

Fix $\varepsilon>0$ and set $Y_1:=B([a,p],\varepsilon)\subset A\setminus NW(C(f))$ and $Y_2:=B([q,a],\varepsilon)\subset A\setminus NW(C(f))$. Consider the line
$$l:=\{[x,p]\mid d(x,a)<\varepsilon\}\subset Y_1.$$
Notice that
$$C(f)^n(l)=\{[f^n(x),f^n(p)]\mid x\in B(a,\varepsilon)\}.$$
Since $f^n(p)\to a$ and $f^n(B(a,\varepsilon))\to S^1$, when $n\to\infty$, there exists $n_0$ such that for all $n\ge n_0$ both $d(f^n(p),a)<\varepsilon$ and $f^n(B(a,\varepsilon))\ni q$ hold. Denote by $x_1\in B(a,\varepsilon)$ such that $f^n(x_1)=q$.
We conclude that $d_H([f^n(p),f^n(x_1)],[a,q])<\varepsilon$, therefore $C(f)^n(Y_1)\cap Y_2\neq\emptyset$.

The dynamics of $C(f)|_B$ is the following:
\begin{itemize}
 \item $y\neq a$ $\Rightarrow$ $C(f)^n([x,y])\stackrel{n\to\pm\infty}{\longrightarrow}S^1$
 \item $[x,y]=\{a\}$ $\Rightarrow$ $C(f)^n([x,y])\stackrel{n\to\pm\infty}{\longrightarrow}\{a\}$
 \item $y=a$, $x\neq a$ $\Rightarrow$ $C(f)^n([x,y])\stackrel{n\to\infty}{\longrightarrow}\{a\}$ and $C(f)^n([x,y])\stackrel{n\to-\infty}{\longrightarrow}S^1$
 \item $[x,y]=S^1$ $\Rightarrow$ $C(f)^n([x,y])\stackrel{n\to\pm\infty}{\longrightarrow}S^1$.
 \end{itemize}

The same reasoning applies to $C(f)|_B$, so the proof of Case (1) is done.

\noindent{\bf Case (2).} Suppose $Fix(f)=\{a_1,\ldots,a_m\}$ and there are no fixed points between the points $a_j$ and $a_{j+1}$ . Denote by $C_j:=[a_j,a_{j+1}]$. It is obvious that the sets $C_j$ are $f$-invariant, therefore the sets
$$D_{ij}:=\{[x,y]\mid x\in C_i, y\in C_j\}$$ are $C(f)$-invariant. It is also easy to see that all $D_{ij}$ are closed as well as $C(S^1)=\bigcup_{i,j}D_{i,j}$. The proof of \ref{lem:B} is completed if we prove that $h_{pol}(f|_{D_{ij}})=2$ for all $i,j$.

We see that $D_{i,j}$ can be identified with $[a_i,a_{i+1}]\times[a_j,a_{j+1}]$ and
since $f([x,y])=[f(x),f(y)]$, $C(f)$ is conjugated to $f\times f$. Therefore we reduce the problem to the computation of the polynomial entropy of
$$C(f)=f\times f:[a_i,a_{i+1}]\times[a_j,a_{j+1}]\to [a_i,a_{i+1}]\times[a_j,a_{j+1}],\quad f\times f(x,y)=(f(x),f(y)).$$
Since $h_{pol}(f\times g)=h_{pol}(f)+h_{pol}(g)$, we have $h_{pol}(C(f)|_{D_{ij}})=h_{pol}(f|_{a_i,a_{i+1}]})+h_{pol}(f|_{[a_j,a_{j+1}]})$. Similarly as in the proof of \ref{lem:A}, we have $h_{pol}(f|_{[a_i,a_{i+1}]})=h_{pol}(f|_{[a_j,a_{j+1}]})=1$.

\noindent{\bf Case (3).} Suppose $a:=a_1$ and $b:=a_2$ are the only two fixed points. Then either $f$ maps both $[a,b]$ and $[b,a]$ to themselves, or one to another. If the latter is the case, then $f^2$ maps both arcs to itself, so we can assume that this is true, and apply the same argument as in (2).

Now we prove the statement for $f^{*k}$. Recall first that $h_{pol}(f)=1$. This can be proved using the coding methods (it is easy to see that $f$ possesses no two mutually singular points), or, alternatively, by refering to Theorem 2 in~\cite{L3}, which states that polynomial entropy of a circle homeomorphism $f$ is $1$ if and only if $f$ is not conjugated to a rotation.

Define a relation $\le$ on $S^1$ by identifying $S^1$ with $[0,1)$. We can again assume that $f$ is orientation preserving, since if not, $f^2$ is and
$$h_{pol}\left(\left(f^2)^{*k}\right)\right)=h_{pol}\left(\left(f^{*k}\right)^2\right)=h_{pol}\left(f^{*k}\right).$$ So we can consider $f$ as an increasing homeomorphism of $[0,1]$ with $f(0)=0$ and $f(1)=1$ (we can assume that $0$ is a fixed point of $f$). As before, by considering a semi-conjugacy $\pi:\left(S^1\right)^k\to\left(S^1\right)^{*k}$, we derive $h_{pol}\left(f^{*k}\right)\le k$. The case $k=2$ is also covered in the same way as for an interval.

If $f$ possesses at least two fixed points, $a$ and $b$, we can define $B$ as a subset of $\left(S^1\right)^{*k}$ consisting of sets of points from the interval $[a,b]$ and identify the map $f^{*k}|_B$ with $\left(f|_{[a,b]}\right)^{*k}$. Therefore we have
$$h_{pol}\left(f^{*k}\right)\ge h_{pol}\left(f^{*k}|_B\right)=h_{pol}\left(\left(f|_{[a,b]}\right)^{*k}\right)=k$$
(the last equality follows from the proof of \ref{lem:A}), so the proof is finished.

If $f$ has only one fixed point, $0$, we can define the sets $\tilde{A}(k)$ and $\hat{A}(k)$ as in the proof of \ref{lem:A}, conclude that $h_{pol}\left(f^k;\tilde{A}(k)\right)=k$, and then prove that $h_{pol}\left(f^{*k};\hat{A}(k)\right)=k$, as in the proof of Lemma~\ref{lem}.

The last statement, $h_{pol}\left(2^f\right)=\infty$, follows in the same way as in the case of an interval. \qed

\begin{rem} If $(X,f)$ and $(Y,g)$ are two dynamical systems and there exists a semi-conjugacy $\pi:X\to Y$ that is uniformly finite-to-one (meaning that there exists $C$ such for any $y\in Y$ it holds $\sharp\;\left(\pi^{-1}(y)\right)\le C$), then the topological entropy $h_{top}(f)$ and $h_{top}(g)$ coincide (see, for example~\cite{R}). It easily follows from this that $h_{top}(f^{*k})=h_{top}(f^{\times k})$ (see~\cite{KO}). An analogous formula for the polynomial entropy is still not proved or disproved, so we had to use the inequality relation, as a particularity of one-dimensional sets.
\end{rem}

\end{document}